\documentclass{ws-ijbc}
\begin{document}

\markboth{}{}

\title{Chaos suppression in a Gompertz-like discrete system of fractional order}
\author{MARIUS-F. DANCA}

\address{Dept. of Mathematics and Computer Science, Avram Iancu University of Cluj-Napoca, Romania\\Romanian Institute od Science and Technology,\\
400487 Cluj-Napoca, Romania,\\
danca@rist.ro}

\author{MICHAL FE\v{C}KAN}
\address{Dept. of Mathematical Analysis and Numerical Mathematics, Faculty of Mathematics, Physics and Informatics, Comenius University in Bratislava, Slovak Republic,\\
Mathematical Institute, Slovak Academy of Sciences, Slovak Republic\\
Michal.Feckan@fmph.uniba.sk}

\maketitle

\author{}

\maketitle

\begin{abstract}
In this paper we introduce the fractional-order variant of a Gompertz-like discrete system. The chaotic behavior is suppressed with an impulsive control algorithm. The numerical integration and the Lyapunov exponent are obtained by means of the discrete fractional calculus. To verify numerically the obtained results, beside the Lyapunov exponent, the tools offered by the 0-1 test are used.
\end{abstract}

\emph{Keywords: Discrete Fractional Order System; Gompertz-like system; Chaos suppression}
\bigskip

\section{Introduction}
In the 19th century Benjamin Gompertz originally derived the Gompertz curve to estimate human mortality \cite{g1}. Later, the Gompertz equation was used to describe growth processes \cite{g2}. One of the re-parametrisation of the Gompertz model is \cite{g3}

\begin{equation*}
\frac{1}{N(t)}\frac{dN(t)}{dt}=a \ln\Big(\frac{N_0}{N(t)}\Big),
\end{equation*}
where $N(t)$ represents the number of cells of the tumor.

The disadvantage of the Gompertz model as remarked by Ahmed in \cite{g4} is that it is not biologically motivated. Therefore, Ahmed introduced the following dynamical equation:
\begin{equation}\label{unu}
\frac{dN(t)}{dt}=-a N(t)+b N(t)^{2/3}.
\end{equation}
For an adequate choice of parameters $a$ and $b$, the growth described by \eqref{unu} is very close to the Gompertz equation in the observable region.

The discrete variant of \eqref{unu} is
\begin{equation}
N_{n+1}=(1-a)N_n+bN_n^{2/3}.
\end{equation}
which, after rescaling, becomes \cite{g4}
\begin{equation}\label{doi}
N_{n+1}=m(N_n^{2/3}-N_n),
\end{equation}

Next, in \cite{g4} it is shown numerically that \eqref{doi}, via a period doubling bifurcation scenario, presents a chaotic behavior.

On the other side, in \cite{dan1} the equation \eqref{doi} is scaled to obtain the following discrete system in a more numerically accessible form:
\begin{equation}\label{trei}
x_{n+1}=6.75r(x_n^\frac{2}{3}-x_n).
\end{equation}

Let next
\begin{equation}\label{ecus}
f_r(x)=6.75r(x^\frac{2}{3}-x).
\end{equation}
 Like the logistic map, on $[0,1]$, $f_r$ is unimodal (one of the main chaos ingredients), because $f_r(0)=f_r(1)=0$, $f_r$ admits only one extreme (maximum), located inside $[0,1]$, at the critical point $x=8/27$, $f_r$ is strictly increasing on $[0,8/27]$ and strictly decreasing on $[8/27,1]$. However, unlike the logistic map, $f_r$ has no negative Schwarzian derivative \cite{schr}.

The chaotic behavior of the system \eqref{trei} can be suppressed using the following periodic variable perturbations \cite{dan1}:

\begin{equation}\label{eee1}
x(n+1)=\begin{cases} f(x(n)), &n\in \mathbb{N},\\
                   (1+\gamma)x(n+1) &\textrm{if } \mod(n,\delta)=0,
     \end{cases}
\end{equation}
where $f$ is some map defining a discrete system ($f_r$ for the considered system), $\gamma\in\mathbb{R}$ represents a relative small perturbation which is applied at every $\delta$ steps. As shown in \cite{dan1}, for adequate $\gamma$ and $\delta$, the chaotic behavior can be suppressed.

This simple control algorithm, introduced first in \cite{gm}, has been adapted for continuous systems of fractional order (FO) \cite{gm2}, discontinuous systems of integer order (IO) \cite{gm5}, discontinuous systems of FO \cite{gm3,gm4}.

In this paper the algorithm \eqref{eee1} is applied to the FO variant of \eqref{trei} to suppress the chaotic behavior. The numerical tools utilized to verify the results are: bifurcation diagram, time series, Lyapunov exponents and the 0-1 test.

The paper is organized as follows: in Section 2 is presented the discrete FO variant of the system \eqref{trei} and in Section 3 presents the numerical results. Section 4 closes the paper and the Appendix section describes briefly the 0-1 test.

\section{The discrete FO system \eqref{trei}}
In order to introduce the FO discrete variant of the system \eqref{trei}, let $q\in(0,1)$, $\mathbb{N}_{1-q}=\{1-q,2-q,3-q,\cdots\}$, $0<q\leq1$ and $f\in C(\mathbb{R},\mathbb{R})$ some discrete map. Therefore, FO systems are modeled by the following initial value problem:
\begin{equation}\label{e1}
\triangle_*^q x(k)=f(x(k-1+q)),\quad k\in \mathbb{N}_{1-q}, ~~x(0)=x_0,
\end{equation}
where $\triangle_*^q v(k)$ is the Caputo delta fractional difference \cite{bal2,14,15}.

With $f:=f_r$ the initial value problem \eqref{e1} for the considered system, \eqref{trei}, becomes
\begin{equation}\label{ex}
\triangle_*^q x(k)=6.75r\big(x(k+q-1)^\frac{2}{3}-x(k+q-1)\big)\quad k\in \mathbb{N}_{1-q},~~x(0)=x_0.
\end{equation}

Hereafter, the FO discrete equations are considered with initial condition $x(0)=x_0$.

The discrete integral form of \eqref{e1} as presented in \cite{bal2} is
$$
x(n)=x(0)+\frac{1}{\Gamma(q)}\sum_{j=1-q}^{n-q}\frac{\Gamma(n-j)}{\Gamma(n-j-q)}f(x(j-1+q)),
$$
or, similarly

\begin{equation*}%\label{e2}
x(n)=x(0)+\frac{1}{\Gamma(q)}\sum_{j=1}^{n}\frac{\Gamma(n-j+q)}{\Gamma(n-j+1)}f(x(j-1)),~~n=1,2,...
\end{equation*}
Therefore, the numerical integration of \eqref{ex} yields
\begin{equation}\label{ex2}
x(n)=x(0)+\frac{6.75r}{\Gamma(q)}\sum_{j=1}^{n}\frac{\Gamma(n-j+q)}{\Gamma(n-j+1)}\big(x(j-1)^\frac{2}{3}-x(j-1)\big),~~n=1,2,...,
\end{equation}
which will be the used as the numerical variant of the FO system \eqref{ex}.

\begin{remark}\label{rem}
Note that due to the discrete memory effect (the present status depends on the all previous
information), one of the main impediments to implement numerically \eqref{ex2}, is the divergency of the term $\sum_{j=1}^n\Gamma(n-j+q)/\Gamma(n-j+1)$.
Thus, it can be proved that
\begin{equation}\label{rela}
\bigg | \sum_{j=1}^{n}\frac{\Gamma(n-j+q)}{\Gamma(n-j+1)}-\frac{n^q}{q}\bigg |\leq\frac{1}{q},~~~ n=1,2,...
\end{equation}
which means that $\sum_{j=1}^{n}\frac{\Gamma(n-j+q)}{\Gamma(n-j+1)}$ and $(n^q)$ grow similarly.
Therefore, for relatively large values of $n$, small errors in the steps of \eqref{ex2} may lead to final computational errors. A simple ``trick'' to increase computationally the iteration number, from only few hundreds to much larger values, is to replace $\Gamma(n-j+q)/\Gamma(n-j+1)$ with $e^{\ln(\Gamma(n-j+q))-ln(\Gamma(n-j+1))}$.
\end{remark}

\begin{figure}
\begin{center}
\includegraphics[scale=1]{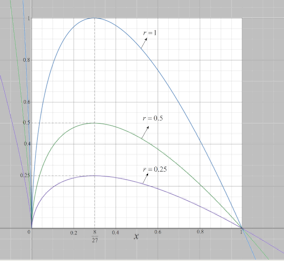}
\caption{Graph of $f_r:[0,1]\rightarrow[0,1], f_r=6.75r(x^{\frac{2}{3}}-x)$, for different values of $r$.}
\label{fig1}
\end{center}
\end{figure}

One can see that considering $f_r$ defined on $[0,1]$, $f_r:[0,1]\rightarrow[0,1]$, and $r\in[0,1]$ the bifurcation parameter, the map graph is like for logistic map within $x\in[0,1]$ (see the graphs for few values of $r$ in Fig. \ref{fig1}). However, as revealed by the numerical results, $x_n$ given by \eqref{ex2}, exceeds the value $x=1$, namely $x_n\in[0,1.5]$.

Due the discrete memory effect of the numerical approach, the Jacobian matrix necessary for the finite-time local Lyapunov exponent, denoted hereafter LE, cannot be obtain directly as for the IO systems. However, by using the natural linearization of \eqref{ex2} along the orbit $x_n$ \cite{bal3}, one has

\begin{equation}\label{d2}
a(n)=a(0)+\frac{6.75r}{\Gamma(q)}\sum_{j=1}^{n}\frac{\Gamma(n-j+q)}{\Gamma(n-j+1)}a(j-1)\bigg(\frac{2}{3}x(j)^{-\frac{1}{3}}-1\bigg),~~a(0)=1.
\end{equation}
Next, the numerical determined LE, $\lambda$, is obtained as follows:
\begin{equation*}
\lambda(x_0)\simeq\frac{1}{n}\ln|a(n-1)|.
\end{equation*}

\begin{figure}
\begin{center}
\includegraphics[scale=1]{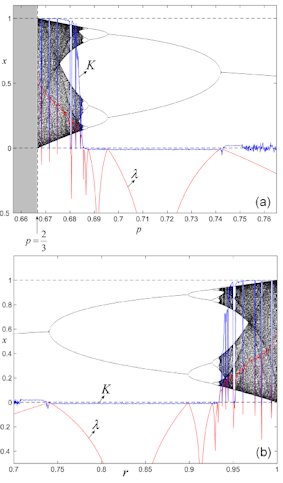}
\caption{Bifurcation diagrams of the map $f_{r,p}$. (a) Bifurcation diagram of the $f_{r,p}$ versus $p$; (b) Bifurcation diagram of $f_{r,p}$ versus $r$.}
\label{fig2}
\end{center}
\end{figure}

Let $f_{r,p}(x)=6.75r(x^p-x)$, as depending on two parameters $p,r$, with $p\in[0.66,0.765]$. To obtain a visual perspective of the influence of the power exponent $p$ and $r$ in the IO system \eqref{trei}, consider the bifurcation diagrams versus $p$ and also versus $r$ (Fig. \ref{fig2}). The LE (red plot) and $K$ value (blue plot) given by the 0-1 test (Appendix) are superimposed on all bifurcation diagrams. As known, if a system behaves regularly, $K$ is approximatively zero, while in the case of chaotic dynamics, $K$ is approximatively 1.

The grey fill in bifurcation diagrams indicates the ranges where the system orbits diverge.

Together with the positive values of the LE, the values of $K$ close to $1$ indicate the chaotic windows, while the negative values of LE and $K$ show the numerical periodic windows (see next section for more details). Note that the values of $K$, especially for the FO system, present some relative small errors (or order $1e-2$ for the 0 value and $1e-3$ for $1$ value), actually typical for the 0-1 test.

\begin{figure}
\begin{center}
\includegraphics[scale=.8]{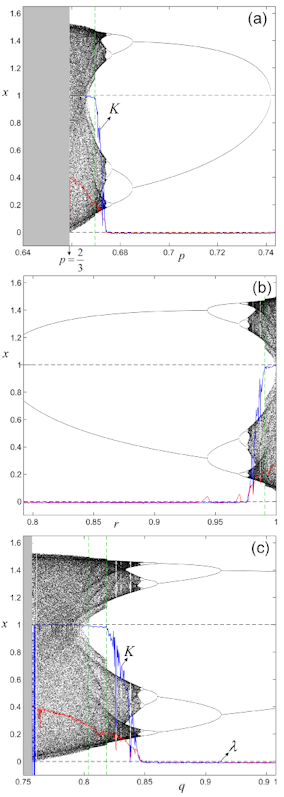}
\caption{Bifurcation diagrams of the FO system \eqref{ex}; (a) Bifurcation diagram versus $p$; (b) Bifurcation diagram versus $r$; (c) Bifurcation diagram versus $q$.}
\label{fig3}
\end{center}
\end{figure}

As Figs. \ref{fig2} a,b and Figs. \ref{fig3} a,b show, there is a significant difference between the values of the LE of the IO and FO (compare, e.g., with Fig. 2 in \cite{bal3}). This difference could be probably explained by the different dynamics in the FO case compared with the IO case, but also by the memory history embedded by the relation \eqref{d2} (see also Remark \ref{rem}).

A common characteristic related to $p$, is the fact that both IO and FO systems present reverse period doubling bifurcations versus $p$, indicating the chaos extinction with the increase of $p$. In the FO case one can see exterior crises (green vertical dotted lines in Fig. \ref{fig3}).

The bifurcation diagram versus $p$, unveils the fact that, for the genuine IO system \eqref{trei}, $p=2/3$ has been chosen as the smallest admissible value (for which chaos is strongest). In this paper, $p$ will be set $p=2/3$ (Fig. \ref{fig3} (a)).

\section{Chaos suppressing}
Using the common notation $x_n=x(n)$, by combining the chaos control algorithm \eqref{eee1} with the numerical integral \eqref{ex2}, one obtains the following algorithm:

\begin{equation}\label{xe1}
x_{n+1}=\begin{cases} x_0+\frac{6.75r}{\Gamma(q)}\sum_{j=1}^{n}\frac{\Gamma(n-j+q)}{\Gamma(n-j+1)}\big(x_{j-1}^\frac{2}{3}-x_{j-1})\big), &n\in \mathbb{N}\\
                   (1+\gamma)x_{n+1}, &\textrm{if } \mod(n,\delta)=0.
     \end{cases}
\end{equation}
This means that at every $\delta$ step (i.e. $n$ is multiple of $\delta$), $x_{n+1}$ is perturbed with $1+\gamma$.

\begin{remark}\label{rufu}
Like continuous-time systems of FO, discrete systems of FO cannot have any nonconstant periodic solution \cite{mich}. Therefore one cannot consider that the system \eqref{ex} admits stable cycles, but only \emph{numerically stable periodic orbits} (NSPO) \cite{dada}, i.e. closed trajectories in the phase space in the sense that the closing error is within a given bound of $1E-n$, with $n$ a sufficiently large positive integer.
\end{remark}

To implement numerically the algorithm \eqref{xe1}, suppose the parameter $r$ in the FO discrete system \eqref{ex} is set such that the system evolve chaotically (here, $r=1$). Then, by choosing adequate values of $\gamma$ and $\delta$, the algorithm can suppress the chaotic behavior, forcing the system to evolve along an NSPO.
Thus, fixing $\delta$ in \eqref{xe1} to some value, to obtain the algorithm parameters values $\gamma$ which suppress the chaos, one determines the bifurcation diagram versus $\gamma\in[\gamma_1,\gamma_2]$ where $\gamma_{1,2}$ have some relative small values. In this paper $\gamma_1=-0.1$ and $\gamma_2=0.1$\footnote{The algorithm can be tested empirically too, by testing values of $\gamma$ or $\delta$ until the chaos is suppressed.}. All experiments have been realized for $q=0.8$.

The moment $n^*$ when the algorithm is applied, is marked with vertical dotted lines in the graphs of time series. $1000$ iterations have been considered.

To identify numerically the regular dynamics obtained with the algorithm (see Remark \ref{rufu}), beside the LE and time series, the data given by the 0-1 test are utilized. Thus, the ranges of $\gamma$ where the LE is not positive, $K$ is (close to) zero, the graphs of $p$ and $q$ are disc-like and the graph of $M$ is not divergent (see Appendix), represent the admissible values of $\gamma$ for chaos suppressing. Beside the LE (red plot), $K$ (blue plot), in the time series, the elements of the NSPOs are plotted red.

\begin{figure}
\begin{center}
\includegraphics[scale=1]{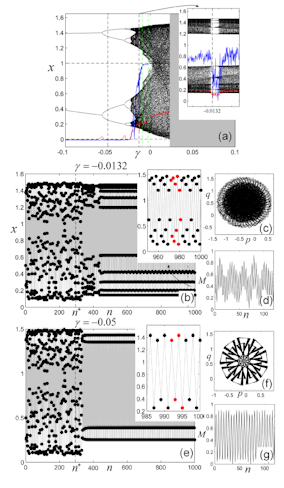}
\caption{Chaos suppression of the FO system \eqref{ex} for $\delta=1$, applied after $n=n^*$; (a) Bifurcation diagram versus $\gamma$. The zoomed area in the time series reveals the chosen value of $\gamma$: $\gamma=-0.0132$; (b) Time series revealing the ten elements of the NSPO (red plot); (c) Graph of $p$, $q$; (d) Graph of $M$; (e) Chaos suppression for $\gamma=-0.05$. Time series reveals the four elements of the NSPO (red plot); (f) Graph of $p$, $q$; (g) Graph of $M$.}
\label{fig4}
\end{center}
\end{figure}

Consider first $\delta=1$, i.e., at every step, $x_{n+1}$ is perturbed with $1+\gamma$.
The bifurcation diagram versus $\gamma$ (Fig. \ref{fig4} (a)) indicates that there are large ranges of $\gamma$ (numerical periodic windows) for which, any value $\gamma$ chosen there, generates an NSPO. For example, for $\gamma=-0.0132$ (see the zoomed area), the system is forced to evolve along the NSPO of 10-period (see the red plot in zoomed area of the time series in Fig. \ref{fig4} (b)). As can be seen in Fig. \ref{fig4} (a), for $\gamma$ within a small neighborhood of $-0.0132$, LE and $K$ are close to zero ($K=0.0033$), the graph of $p$ and $q$ is disk-like and the mean-square displacement $M$ is not unbounded, underlying the numerically periodic motion. More clear values of LE and $K$ can be obtained for $\gamma$ chosen in the large numerical periodic window, $\gamma=-0.05$ (Fig. \ref{fig4} (c)). Now $K=-0.0035$. This happens probably due to an inertia-like phenomenon necessary to LE and especially to $K$ to stabilize for the considered parameter ranges.

Note that at least two crisis can be remarked (vertical green dotted lines).

\begin{figure}
\begin{center}
\includegraphics[scale=1]{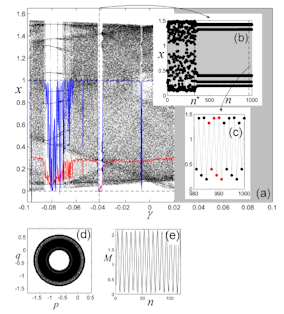}
\caption{Chaos suppression of the FO system \eqref{ex} for $\delta=3$ and $\gamma=-0.04$; (a) Bifurcation diagram versus $\gamma$; (b) Time series for $\gamma=-0.04$; (c) Zoomed area revealing the six elements of the NSPO (red plot); (d) Graph of $p$, $q$; (e) Graph of $M$.}
\label{fig5}
\end{center}
\end{figure}

For e.g. $\delta=3$, i.e., $x_{n+1}$ is perturbed only every third steps, there still exist numerically periodic windows where the chaos is suppressed but, as expected, their size is smaller (see Fig. \ref{fig5}, where for $\gamma=-0.04$ a numerically 5-period orbit is obtained) and $K=-0.0048$. Again, the graph of $p$ and $q$ is disk-like and $M$ has a bounded evolution.

\begin{figure}
\begin{center}
\includegraphics[scale=1]{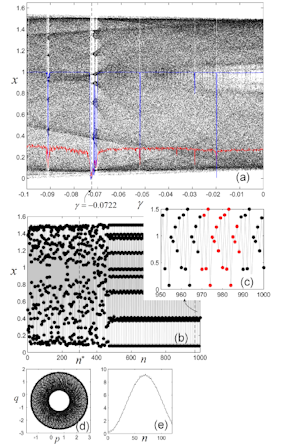}
\caption{Chaos suppression of the FO system \eqref{ex} for $\delta=5$ and $\gamma=-0.0722$; (a) Bifurcation diagram versus $\gamma$; (b) Time series with a zoomed area showing the 19-period orbit of the obtained NSPO (red plot in Fig. \ref{fig6} (c)); (d) Graph of $p$, $q$; (e) Graph of $M$.}
\label{fig6}
\end{center}
\end{figure}

The largest value of $\delta$ for which the system still can be controlled is $\delta=5$ (Fig. \ref{fig6}). For $\gamma=-0.0722$ an NSPO with a large 19-period is obtained. In this case $K=0.0317$.

For larger values of $\delta$, no NSPOs are found.

\begin{figure}
\begin{center}
\includegraphics[scale=.7]{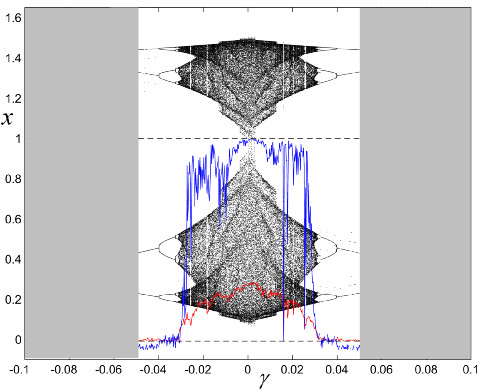}
\caption{Chaos control of the FO system \eqref{ex} obtained with the algorithm \eqref{pluss}, $\delta=2$. The bifurcation diagram versus $\gamma$ shows a symmetry allowing beside negative for $\gamma$, positive values too.}
\label{fig7}
\end{center}
\end{figure}

Note that other possible variants of \eqref{xe1} can be considered and implemented in practice, such as
\begin{equation}\label{pluss}
x_{n+1}=\begin{cases} f(x_n), &n\in \mathbb{N},\\
                   x_{n+1}+\gamma &\textrm{if } \mod(n,\delta)=0,
     \end{cases}
\end{equation}
which, for e.g. $\delta=2$ and $|\gamma|\leq0.5$, has the symmetric bifurcation diagram in Fig. \ref{fig7}.

\section{Conclusion and discussion}
In this paper we introduced the FO variant of a biological discrete system, modeled with Caputo's derivative. Also, a control algorithm to suppress the chaotic behavior is presented. The numerical integration of the system was made by using the discrete integral presented in \cite{bal2}. Under the action of the control algorithm, the chaotic behavior of the system can be transformed into regular motion, i.e. orbits which are numerically apparent stable periodic (as known, both continuous and discrete FO systems admit no stable periodic orbits).

Every $\delta$ steps, the algorithm impulses periodically the system variable $x_{n+1}$, in such a way that $x_{n+1}=x_{n+1}(1+\gamma)$ with $\gamma$ chosen from the bifurcation diagram versus $\gamma$.

Note that this kind of perturbations can rarely happen in nature, but with positive probability, when some system accidently perturbs periodically its variables.

While for most of continuous-time systems, there are relatively small differences between the IO and their FO variants, in the case of the discrete system \eqref{trei} there are some notable differences which could also characterize other discrete systems. Thus, compared to IO variant of the system, \eqref{trei}, for which the LE has an expected evolution for IO discrete systems (with relative large negative and positive values (Fig. \ref{fig2})), in the case of the considered FO system \eqref{ex}, for all considered numerical experiments, the negative values of the LE are actually constant and close to zero for relatively large ranges of $\gamma$ (see also \cite{xin,bal3}). The question if the LE is really negative, or is zero, remains an open problem.

Like for IO systems, the $K$ value given by the 0-1 test, presents oscillations in some critical ranges of $\gamma$. Also, due to the inherent numerical characteristics of the results of the test 0-1 or, probably due to NSPOs (Remark \ref{rufu}), the errors in calculating the 0 value is relatively larger ($1e-2$) than the errors in calculating 1, which is $1e-3$.

All experiments with the algorithm \eqref{xe1} revealed that only negative values of $\gamma$ allow the chaos suppression, suggesting that the underlying system must free energy every $\delta$ steps to evolve along some NSPO. On the other side, the algorithm \eqref{pluss} allows the chaos suppression for both negative and positive values, which means that the system can be stabilized by either lose or gain energy.

Another open problem regards the performances of the numerical integration of FO discrete systems (see Remark \ref{rem}). Finding and using other numerical integrals of \eqref{e1} could improve the obtained numerical results.

\section*{Appendix}

The '0-1'  test has been developed in \cite{got2}, being designed to distinguish chaotic behavior from regular behavior in continuous and discrete systems.
The input is a time series, the test being easy to implement. Note it does not need the system equations.
Let a discrete or continuous-time dynamical system (of IO or FO) and a one-dimensional observable data set be determined from a time series, $\phi(j)$, $j=1,2,...,N$, with $N$ some positive integer. It is proved that the test states that a nonchaotic motion is bounded, while a chaotic dynamic behaves like a Brownian motion \cite{uu}.
To obtain the four elements generated by the test, the asymptotic growth $K$, $p$, $q$, and the mean-square displacement $M$, the following steps are to be determined:

\noindent 1) First, for $c\in[0,2\pi]$, compute the translation variables $p$ and $q$ \cite{got2}:

\[
p(n)=\sum_{j=1}^n\phi(j)\cos(jc),~~~ q(n)=\sum_{j=1}^n\phi(j)\sin(jc),
\]
for $n=1,2,...,N$. However, $c$ can be chosen for example within a narrow interval $[\pi/5,4\pi/5]$ as mentioned in \cite {got1}.

\noindent 2) Next, in order to determine the growths of $p$ and $q$, the mean-square displacement $M$ is determined:
\[
M(n)=\lim_{N\rightarrow \infty}\frac{1}{N}\sum_{j=1}^N[p(j+n)-p(j)]^2+[q(j+n)-q(j)]^2.
\]
where $n\ll N$ (in practice, $n=N/10$ represents a good choice).

\noindent 3) The asymptotic growth rate $K$ is defined as
\[
K=\lim_{n\rightarrow \infty}\log M(n)/\log n.
\]
If the system dynamics is regular (i.e. periodic or quasiperiodic) then $K \approx 0$, otherwise, if the underlying dynamics is chaotic, $K \approx 1$.

Note that, recently, the 0-1 test has been used successfully to identify strange nonchaotioc attractors \cite{gop}.

\begin{figure}
\begin{center}
\includegraphics[scale=1]{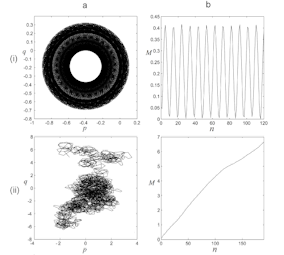}
\caption{The 0-1 test applied to the logistic map $x_{n+1}=rx(1-x)$; (a) Graphs of $p$, $q$; (b) Graphs of $M$; (i) The case $r=3.5$; (ii) The case $r=3.9$.}
\label{fig8}
\end{center}
\end{figure}

Figs. \ref{fig8} present the case of the logistic map $x_{n+1}=rx_n(1-x_n)$. In Figs. (i) and (ii), the cases of $r=3.55$, when the system behaves regularly, and $r=3.9$, when the system behaves chaotically, are considered, respectively. Figs. (a) and (b) show the graph of $p$, $q$ and $M$, respectively.

\begin{thebibliography}{99}                                                                                               %

\bibitem[Nicol et al. (2001)]{uu}
Nicol M., Melbourne I. and Ashwin P. ``Euclidean extensions of dynamical systems'',
  \emph{Nonlinearity} 14~(2) (2001) 275--300.

\bibitem[Gompertz(1825)]{g1} Gompertz, B. [1825] ``On the nature of the function expressive of the law of human mortality, and on a new mode of determining the value of life contingencies'', \emph{Phil. Trans. Roy. Soc. London} \textbf{115}, 513-583.
\bibitem[Winsor(1932)]{g2} Winsor, C. [1932] ``The Gompertx curve as a growth equation'', \emph{Proc. Nat. Acad. Sciences}, \textbf{18}, 1-8.
\bibitem[Swan(1990)]{g3} Swan, G.W. [1990] ``Role of optimal control theory in cancer chemotherapy'', \emph{Math. Biosciences }\textbf{101}, 237-284
\bibitem[Ahmed(1992)]{g4} Ahmed E. [1992] ``Fractals and chaos in cancer models'', \emph{Int. J. Theor. Phys. }\textbf{32}(2), 353-355.
\bibitem[Wu \& Baleanu(2014)]{bal2} Wu, G.-C., Baleanu, D. 2014 ``Discrete fractional logistic map and its chaos'', \emph{Nonlinear Dynam}. \textbf{75}(1-2), 283--287.
\bibitem[Goodrich \& Peterson(2015)]{14} Goodrich, C., Peterson, A.C. [2015] \emph{Discrete Fractional Calculus} (Springer).
\bibitem[Feckan \& Pospisil(2014)]{15} Feckan, M., Pospisil, M. [2014] ``Note on fractional difference Gronwall inequalities'', \emph{Electron. J. Qual. Theo.} \textbf{44}, 1–18
\bibitem[Codreanu \& Danca(1997)]{dan1} Danca, M.-F. and Codreanu, S. [1997] ``Suppression of chaos in a one-dimensional mapping'', \emph{J. Biol. Phys.} \textbf{23}, 1–9.
\bibitem[Guemez \& Matias(1993)]{gm} Guemez, J. and Matias, M.A. [1993] ``Control of chaos in unidimensional maps'', \emph{Phys. Lett. A} \textbf{181}, 29-32.
\bibitem[Danca et al.(2016)]{gm2} Danca, M.-F., Tang, W., Chen, G. [2016] ``Suppressing chaos in a simplest autonomous memristorbased circuit of fractional order by periodic impulses'', \emph{ Chaos Soliton Fract.} \textbf{84}, 31-40.
\bibitem[Danca \& Garrappa(2015)]{gm3} Danca, M.-F., Garrappa, R. [2015] ``Suppressing chaos in discontinuous systems of fractional order by active control'', \emph{Appl. Math. Comput.} \textbf{257}, 89-102.
\bibitem[Wu \& Baleanu(2015)]{bal3}Wu, G.-C., Baleanu, D. [2015] ``Jacobian matrix algorithm for Lyapunov exponents of the discrete fractional maps'', \emph{Commun. Nonlinear Sci.} \textbf{22}(1-3), 95-100.
\bibitem[Danca(2012)]{gm4} Danca, M.-F. [2012] ``Chaos suppression via periodic change of variables in a class of discontinuous dynamical systems of fractional order'', \emph{Nonlinear Dynam. } \textbf{70} (1), 815-823.
\bibitem[Marek \& Schreiber(1991)]{schr}Marek M. and Schreiber I. [1991] \emph{Chaotic Behavior of Deterministic Dissipative Systems} (Cambridge Univ. Press).
\bibitem[Diblik et al.(2015)]{mich}Diblik, J., Fe ckan, M., Posp\'{\i}\v{s}il, M. ``Nonexistence of periodic solutions and S-asymptotically periodic solutions in fractional diff erence equations'', \emph{Appl. Math. Comput.} \textbf{257}, 230-240.
\bibitem[Danca et al.(2018)]{dada}Danca M.-F., Feckan M., Kuznetsov N. and Chen G. [2018] ``Fractional-order PWC systems without zero Lyapunov exponents'', \emph{Nonlinear Dynam.}, \textbf{92}(3), 1061–1078.
\bibitem[Gottwald \& Melbourne(2004)]{got2} Gottwald G., Melbourne I. [2004] ``A new test for chaos in deterministic systems'',
  Proceedings of the Royal Society A: Mathematical, Physical and Engineering
  Sciences \textbf{460}(2042)603--611.
\bibitem[Gottwald \& Melbourne(2009)]{got1}
Gottwald G., Melbourne I. [2009] ``On the implementation of the 0-1 test for chaos'',
  \emph{SIAM J. Appl. Dyn. Syst. }\textbf{8}(1) 129--145.
\bibitem[Gopal et al.(2013)]{gop}Gopal R., Venkatesan, A. and Lakshmanan M. [2013] ``Applicability of 0-1 test for strange nonchaotic attractors'', \emph{Chaos} \textbf{23}, 023123.
\bibitem[Xin et al.(2017)]{xin}Baogui Xin, Li Liu, Guisheng Hou, Yuan Ma [2017] ``Chaos Synchronization of Nonlinear Fractional Discrete Dynamical Systems via Linear Control'', \emph{Entropy} \textbf{19}(7), 351.
\bibitem[Danca(2004)]{gm5}Danca M.-F. [2004] ``Controlling chaos in discontinuous dynamical systems'', \emph{Chaos Soliton Fract.} \textbf{22}(3) 605-612.
\end{thebibliography}
\end{document}